\documentclass[12pt]{amsart}
\usepackage{amsmath,amsthm,amssymb}
\usepackage{enumitem}
\usepackage{multicol, hyperref}
\usepackage{amsfonts}
\usepackage[utf8]{inputenc}
\usepackage{graphicx,graphics}
\usepackage{color,mathrsfs}

\newtheorem{thm}{Theorem}
\newtheorem*{thm*}{Theorem}

\theoremstyle{definition}

\newcommand{\Q}{\mathbb{Q}}
\newcommand{\Qp}{\mathbb{Q}_p}
\newcommand{\R}{\mathbb{R}}
\newcommand{\CC}{\mathbb{C}}

\newcommand{\p}{\mathfrak{p}}

\begin{document}
\title{Explicit formula in a imaginary quadratic number field}
\author{O. F. Casas-S\'anchez}
\address{Escuela de Matem\'aticas y Estad\'{\i}stica, UPTC, Tunja}
\email{oscar.casas01@uptc.edu.co}
\author{J. J. Rodr\'{\i}guez-Vega}
\address{Departmento de Matem\'aticas, Universidad Nacional de Colombia, Bogot\'a}
\email{jjrodriguezv@unal.edu.co}

\begin{abstract}
In \cite{Haran_1990}, Haran, using Riesz potentials, presents a version of the classical explicit formula for the Riemann zeta function that treats all places equally. In this article, we extend Haran's results to the case of an imaginary quadratic extension of $\Q$. To this end, we define Riesz kernels for both the totally ramified and complex cases. Although Haran implicitly addressed the unramified extension, we also include this case for the sake of completeness. Following Haran's approach with the Riemann zeta function, we demonstrate that these Riesz kernels naturally arise in connection with the contributions of the totally ramified and complex places to the Dedekind zeta function of the imaginary quadratic extension.
\end{abstract}\maketitle

\section{Introduction}
Let $\zeta(s)=\prod_{p\leq \infty}\zeta_p(s)$ denote the classical Riemann zeta function. Riemann knew the exact relation between the zeros of $\zeta(s)$ and the distribution of primes, which is known as the explicit formula. However, it was Weil who shed new light on this with the Weil explicit formula \cite{Weil_1952}, \cite{weil_1972}. Given a smooth and compactly supported function $f:\R^+\rightarrow\R$, let $\mathcal{M}^s(f)$ represent its Mellin transform. By the functional equation $\zeta(s)=\zeta(1-s)$, we have:\begin{align*}
	\sum_{\zeta(s)=0} &\mathcal{M}^s(f)-\mathcal{M}^0(f)-\mathcal{M}^1(f)=\frac{1}{2 \pi i} \oint \hat{f}(s) d \log \zeta(s)\\
	&= \sum_{p \leq \infty} \frac{1}{2 \pi i} \int_{\frac{1}{2}-i \infty}^{\frac{1}{2}+i \infty} \hat{f}(s) d \log \frac{\zeta_p(s)}{\zeta_p(1-s)} \stackrel{\text { def }}{=} \sum_{p \leq \infty} W_p(f)\\
\end{align*}
Among many extensions of Weil's work, Haran provided a new interpretation of $W_p(f)$ using Riesz potentials \cite{Haran_1990}. Haran showed that all the local contributions to the explicit formula can be expressed in terms of the infinitesimal generator of the Riesz kernel associated with the corresponding local field.

In this article, we extend Haran's results to the case of an imaginary quadratic number field. To achieve this objective, we define the Riesz potential for the complex field $\CC$ and for the corresponding extensions of the local fields $\Q_p$, which include unramified and totally ramified quadratic extensions. Then, we demonstrate that we can glue together this operator in an operator over the ideles and show the corresponding explicit formula.

It is worth noting that Burnol also proposed an interpretation of the Weil explicit formula \cite{Burnol_2000}.

\subsection{Main results}
Let $K$ be an imaginary quadratic number field, and denote by $\mathcal{O}_K$ its ring of integers, $D_K$ its discriminant, and $N = N_{K/\Q}$ the norm map. The complete Dedekind zeta function of $K$ is given by
\begin{equation*}
	\zeta_K(s)=\prod_{v}\zeta_v(s),
\end{equation*}
the product taken over all places $v$ of $K$. The local factors are given as follows: If $v=\p$ a prime ideal of $K$, then
\begin{equation*}
	\zeta_v(s)=(1-N(\p)^{-s})^{-1},
\end{equation*}
if $v$ is the real place, then
\begin{equation*}
	\zeta_v(s)=(\pi)^{s/2}\Gamma\bigl(\frac{s}{2}\bigr).
\end{equation*}
if $v$ is the complex place, then
\begin{equation*}
	\zeta_v(s)=(2\pi)^{1-s}\Gamma(s).
\end{equation*}

$\zeta_K(s)$ satisfies the functional equation
\begin{equation*}
	\zeta_K(s)=|D_K|^{\frac{1}{2}-s} \zeta_K(1-s).
\end{equation*}

We denote by $K_v$ the completion of $K$ at the place $v$. In the case of a finite place $v=\p$ above the rational prime $p$ ($p\neq 2$), we have the following posibilities depending of the Jacobi symbol $\bigl(\frac{D_K}{p}\bigr)$. Recall that $\Qp$ is the field of $p-$adic numbers:
\begin{itemize}
\item If $\bigl(\frac{D_K}{p}\bigr)=0$, then $p$ ramifies, $N(\p)=p$, and $K_v/\Q_p$ is a totally ramified quadratic extension.
\item If $\bigl(\frac{D_K}{p}\bigr)=1$, then $p$ splits, $p=\p\overline{\p}$, $N(\p)=N(\overline{\p})=p$, and $K_v=\Q_p$.
\item If $\bigl(\frac{D_K}{p}\bigr)=-1$, then $p$ is inert, $N(\p)=p^2$, and $K_v/\Q_p$ is an unramified quadratic extension. In fact, $K_{v}=\Q_p(\sqrt{\tau})$, where $D_K\equiv \tau \pmod{p}$,
\end{itemize}
and in the archimedean case we denote $K_{\infty}=\CC$.

We choose an absolute value $|\cdot |_v$ on $K_v$, normalized by $|\p |_v=N(\p)^{-1}$ if $v$ is non-archimedean, and by the square of the standard absolute value in the case $K_\infty$. Let $\mathcal{O}_v=\{x\in K_v\mid |x|_v\leq 1\}$ be the ring of integers of $K_v$, and let $\mathcal{O}_v^{\times}=\{x\in K_v\mid |x|_v=1\}$ be the group of units. We denote by $\psi_v$ the canonical character of $K_v$, given by $\psi_{\infty}(z)=e^{-2\pi z\overline{z}}$ for $v=\infty$, and by $\psi_v(x)=\psi_p(Tr_{K_v/\Q_p}(x))$ for $v=\p$ above $p$, where $Tr_{K_v/\Q_p}$ is the trace map, and $\psi_p$ is the canonical character of $\Q_p$ given by $\psi_p: \mathbb{Q}_p \rightarrow \mathbb{C}^{\times}$. $\psi_p(x)=e^{-2\pi i \{x\}_p }$, where $\{x\}_p$ is the fractional part of $x\in\Q_p$ (see \cite{vladimirov1994})

We use $dx_v$ to denote twice the Lebesgue measure for $v=\infty$, and for $v$ finite, the Haar measure normalized by $\int_{\mathcal{O}_v}\,dx_v=N(\mathfrak{d}_v)^{-1/2}$, where $\mathfrak{d}_v$ is the different ideal of $\mathcal{O}_v$, so that every measure is self-dual. We define the Fourier transform of a function $\varphi$ as $\mathcal{F}\varphi(y)=\int_{K_v}\varphi(x)\overline{\psi_v}(xy)\,dx_v$. Additionally, we define $\phi_v(x)$ as follows: $\phi_{\infty}(x)=e^{-2\pi x\overline{x}}$ for $v=\infty$, and $\phi_v(x)=1_{\mathcal{O}_v}(x)$, the characteristic function of $\mathcal{O}_v$, for $v$ finite. Also let $d^{\times}x_v$ the multiplicative Haar measure, defined by $d^{\times}x_\infty=\frac{dx_\infty}{|x|_\infty}$ and for the finite case as $d^{\times}x_v=\frac{dx_v}{(1-N(\p)^{-1})|x_v|_v}$.

Let $\mathbb{A}_K$ represent the ring of adeles of $K$, and $\mathbb{A}_K^\times$ the ideles. The space $\mathscr{S}(\mathbb{A}^{\times}_K) = \otimes{v} \mathscr{S}(K_v^{\times})$ consists of finite linear combinations of elementary functions $f=\otimes_v f_v$, where $f_v$ is a smooth (locally constant for $v<\infty$) and compactly supported function on $K_v^\times$, with $f_v=\phi_v^\times$ for almost all $v$'s. A function is considered symmetric if $f_v$ is invariant under $\mathcal{O}_K^\times$ for all $v$, let
\begin{equation*}
	\mathscr{S}=\Bigl\{\tilde{f}=\sum_{k\in K^\times} f(kx) \mid f\in \mathscr{S}(K_v^{\times}), f \text{ symmetric}  \Bigr\}
\end{equation*}
this consist of all smooth functions on the positive reals that have compact support, via
\begin{equation*}
	K^\times\backslash\mathbb{A}_K^\times/\prod_v \mathcal{O}_v^\times \cong \R^+.
\end{equation*}

Let $f\in \mathscr{S}$, which, based on the former correspondence, can be interpreted as a smooth and compactly supported function on $\mathbb{R}^+$. We define $f|_v(x_v)=f(|x_v|_v)$ as the associated symmetric function on $K_v$ obtained by restriction.

Following Tate's thesis, we denote by $Z_v(\phi_v,s)$ the local zeta function, which is given by:
\begin{align*}
	Z_v&(\phi_v, s)=\\
    &\int_{K_v^{\times}} \phi_v(x)|x|_v^s\,d^{\times}x_v=
	\begin{cases}
		\frac{N(\mathfrak{d}_v)^{-1/2}}{1-N(\p)^{-s}} & \text{ for }K_v/\Q_p \text{ totally ramified}\\ 
		\frac{1}{1-N(\p)^{-s}} & \text{ for }K_v=\Q_p \\ 
		\frac{1}{1-N(\p)^{-s}} & \text{ for }K_v/\Q_p \text{ unramified}\\ 
		(2\pi)^{1-s}\Gamma(s) & \text{ for }v=\infty.
	\end{cases}
\end{align*} 

We also define $Z(\phi,s)=\prod_v Z_v(\phi_v,s)$, where $\phi(x)=\prod_v\phi_v(x)$ is the corresponding function in the adeles $\mathbb{A}_K$. It is well-known (\cite{cassels}) that:
\begin{equation*}
	Z(\phi,s)=Z(\mathcal{F}(\phi),1-s),
\end{equation*}
this implies the functional equation mentioned earlier for $\zeta_K(s)$.

Our first results concerns the Riesz kernels for the fields $K_v$. We define
\begin{equation*}
	R_v^s(x)=\frac{Z_v(\mathcal{F}(\phi_v),1-s)}{Z_v(\phi_v, s)}|x|_v^{s-1}\,dx_v,
\end{equation*}
It is worth noting that according to Tate's thesis \cite{cassels}, the meromorphic function $Z_v(\mathcal{F}(f_v), 1-s)/Z_v(f_v, s)$ is independent of the choice of $f_v$. This function is referred to as the Gamma factor of $K_v$, and in our case, we have the following:
\begin{equation*}
	R_v^s(x)=
	\begin{cases}
		N(\mathfrak{d}_v)^{1/2-s}\dfrac{\zeta_v(1-s)}{\zeta_v(s)}|x|_v^{s-1}\,dx_v & \text{ for }K_v/\Q_p \text{ totally ramified}\\
		\dfrac{\zeta_v(1-s)}{\zeta_v(s)}|x|_v^{s-1}\,dx_v & \text{ for }K_v=\Q_p \\
		\dfrac{\zeta_v(1-s)}{\zeta_v(s)}|x|_v^{s-1}\,dx_v & \text{ for }K_v/\Q_p \text{ unramified}\\ 
		\dfrac{\zeta_\infty(1-s)}{\zeta_\infty(s)}|x|_\infty^{s-1}\,dx_\infty & \text{ for $v$ complex},
	\end{cases}
\end{equation*} 
We prove that the Riesz kernel has a meromorphic continuation and generates a convolution semigroup acting by convolution. In the unramified case, this is already known (see \cite{Kochubei}, \cite{Rodriguez_2008}). However, in the totally ramified and complex cases, our results are new.

Next, we calculate the infinitesimal generator of this group. More precisely
\begin{thm*} 
Let $K$ be an imaginary quadratic number field, and consider any place $v$. Then, for the corresponding Riesz kernel $R_v^s$ defined in $K_v$, we have that for any $f\in\mathscr{S}$ and $\sigma>0$
	\begin{equation*}
		\left.\frac{\partial}{\partial s}\right|_{s=0} R_v^{-s}\ast f|_v(1)=\frac{1}{2\pi i}\int_{\sigma-i\infty}^{\sigma+i\infty}\mathcal{M}^s(f)\,d\log\Bigl(\frac{Z_v(\phi_v, s)}{Z_v(\mathcal{F}(\phi_v),1-s)} \Bigr)
	\end{equation*}
\end{thm*}

In our main result, we renormalize and combine the above kernels to give us an operator $\Delta_{\mathbb{A}^{\times}}^s$ on the ideles $\mathbb{A}^{\times}$ of $K$. We then demonstrate that
\begin{thm*}
Let $f\in\mathscr{S}$, then
\begin{equation*}
\sum_{k \in K^{\times}} \left.\frac{\partial}{\partial s}\right|_{s=0} \Delta_{\mathbb{A}^{\times}}^s f(k)=\sum_{\zeta_K(s)=0} \mathcal{M}^s(f)-\mathcal{M}^0(f)-\mathcal{M}^1(f)
\end{equation*}
In this expression, the right-hand sum extends over all the zeros of $\zeta_K(s)$.
\end{thm*}

\subsection{Outline of the paper}
In Section \ref{RieszPot}, we introduce the Riesz potentials associated with the places of an imaginary quadratic number field, discussing our results on totally ramified extensions and reviewing known results on unramified extensions. Additionally, we introduce the Riesz kernel for the complex field $\CC$. Across all these cases—totally ramified, unramified, and complex—we explicitly calculate the infinitesimal generators. These kernels are shown to form a convolution semigroup. Finally, in Section \ref{ExplicitFormula}, we combine these infinitesimal generators into an adelic operator, which, as Haran demonstrated for $\Q$, is connected with the Dedekind zeta function of $K$.

It is worth mentioning that our results are, in fact, valid for any quadratic extension of $\mathbb{Q}$. We have chosen to focus on the imaginary extensions because the corresponding Riesz kernels have not been studied. In the case of real extensions, our results rely on the well-established real Riesz kernel.


\section{Riesz Potentials}\label{RieszPot}

In this section, we define the Riesz potentials, show  that they generate a convolution semigroup, and calculate their infinitesimal generator.

\subsection{Totally ramified case} Let $K_v/\Qp$ be totally ramified quadratic, in this case the Riesz kernel is given by
\begin{equation*}
	R_v^s(x)=N(\mathfrak{d}_v)^{\frac{1}{2}-s}\dfrac{1-N(\p)^{-s}}{1-N(\p)^{s-1}}|x|_v^{s-1}\,dx_v,
\end{equation*}
$\text{Re } s>0$, $s\not\equiv 1 \pmod{\frac{2\pi i}{\log N(\p)}}$.

As a distribution, it has a meromorphic continuation to all $s$, which is given by
\begin{align*}
    R_v^s(\varphi) &= N(\mathfrak{d}_v)^{-s}\dfrac{1-N(\p)^{-1}}{1-N(\p)^{s-1}}\varphi(0) \\
    &\quad + N(\mathfrak{d}_v)^{\frac{1}{2}-s}\dfrac{1-N(\p)^{-s}}{1-N(\p)^{s-1}}\Bigl ( \int\limits_{|x_v|_v\leq 1} (\varphi(x)-\varphi(0))|x_v|_v^{s-1}\,dx_v \\
    &\quad + \int\limits_{|x_v|_v > 1} \varphi(x)|x_v|_v^{s-1}\,dx_v \Bigr).
\end{align*}

In particular, for $\text{Re } s>0$:

\begin{align*}
	R_v^s(x)&= N(\mathfrak{d}_v)^{\frac{1}{2}-s}\dfrac{1-N(\p)^{-s}}{1-N(\p)^{s-1}}\int_{K_v}\varphi(x) |x|_v^{s-1}\,dx_v \\
	R_v^{-s}(x)&=N(\mathfrak{d}_v)^{\frac{1}{2}+s}\dfrac{1-N(\p)^{s}}{1-N(\p)^{-s-1}}\int_{K_v}\dfrac{\varphi(x)-\varphi(0)}{|x|_v^{1+s}}\,dx_v
\end{align*}
and for $s=0$, $R_{v}^0(\varphi)=\varphi(0)$, i.e. $R_{v}^0=\delta$.

Now, regard $R_v^s$ as an operator via convolution
\begin{multline*}
    R_v^s \ast \varphi(y)= N(\mathfrak{d}_v)^{-s}\dfrac{1-N(\p)^{-1}}{1-N(\p)^{s-1}}\varphi(y) 
    + N(\mathfrak{d}_v)^{\frac{1}{2}-s}\dfrac{1-N(\p)^{-s}}{1-N(\p)^{s-1}}\times \\
    \Biggl(\, \int\limits_{|x_v|_v\leq 1} (\varphi(y+x)-\varphi(y))|x_v|_v^{s-1}\,dx_v \\
    + \int\limits_{|x_v|_v > 1} \varphi(y+x)|x_v|_v^{s-1}\,dx_v \Biggr)
\end{multline*}
once more, for $\text{Re } s>0$:
\begin{align*}
	R_v^s \ast \varphi(y)&= N(\mathfrak{d}_v)^{\frac{1}{2}-s}\dfrac{1-N(\p)^{-s}}{1-N(\p)^{s-1}}\int_{K_v}\varphi(y+x) |x|_v^{s-1}\,dx_v \\
	R_v^{-s} \ast \varphi(y)&=N(\mathfrak{d}_v)^{\frac{1}{2}+s}\dfrac{1-N(\p)^{s}}{1-N(\p)^{-s-1}}\int_{K_v}\dfrac{\varphi(y+x)-\varphi(y)}{|x|_v^{1+s}}\,dx_v
\end{align*}

For $\varphi$ locally constant and with compact support we have that $R_v^s \ast \varphi(y)$ is also locally constant and that $R_v^s \ast \varphi(y)=O(|x|_v^{\text{Re}(s)-1})$, thus we can form the convolution $R_v^{s'} \ast(R_v^s \ast \varphi)=R_v^{s'+s} \ast \varphi$ for $\text{Re }(s'+s)<1$. (see \cite{Kochubei})

In particular, $R_v^{-s}$ for $\text{Re }s\geq 0$ forms a semigroup of operatos whose infinitesimal generator we will calculate below.

\begin{thm}
Let $f\in\mathscr{S}$ and $\sigma>0$, then
\begin{equation*}
		\left.\frac{\partial}{\partial s}\right|_{s=0} R_v^{-s}\ast f|_v(1)=\frac{1}{2\pi i}\int_{\sigma-i\infty}^{\sigma+i\infty}\mathcal{M}^s(f)\,d\log\Bigl(N(\mathfrak{d}_v)^{s-\frac{1}{2}}\frac{1-N(\p)^{s-1}}{1-N(\p)^{-s}} \Bigr)
	\end{equation*}
\end{thm}

\begin{proof}

\begin{equation*}
\begin{aligned}
    \left.\frac{\partial}{\partial s}\right|_{s=0} R_{v}^{-s} &* f|_v(1) = \left.\frac{\partial}{\partial s}\right|_{s=0} \Biggl\{N(\mathfrak{d}_v)^{s} \dfrac{1-N(\p)^{-1}}{1-N(\p)^{-s-1}}f|_v(1) \\
    &\quad + N(\mathfrak{d}_v)^{s+\frac{1}{2}}\dfrac{1-N(\p)^{s}}{1-N(\p)^{-s-1}}\Biggl(\int\limits_{|x_v|_v\leq 1}\dfrac{f|_v(1-x_v)-f|_v(1)}{|x_v|_v^{1+s}}\,dx_v\\&\quad +\int\limits_{|x_v|_v>1}\dfrac{f|_v(1-x_v)}{|x_v|_v^{1+s}}\,dx_v \Biggr) \Biggr\}\\
    &= \log(N(\mathfrak{d}_v))f|_v(1)-\frac{\log(N(\p))}{1-N(\p)^{-1}}\Biggl(N(\p)^{-1}f|_v(1) + N(\mathfrak{d}_v)^{\frac{1}{2}}\times \\
    &\quad \Biggl(\int\limits_{|x_v|_v\leq 1}\dfrac{f|_v(1-x_v)-f|_v(1)}{|x_v|_v}\,dx_v  +\int\limits_{|x_v|_v>1}\dfrac{f|_v(1-x_v)}{|x_v|_v}\,dx_v\Biggr)\Biggr)\\
    \hspace{1.5cm}&= \log(N(\mathfrak{d}_v))f|_v(1)-\frac{\log(N(\p))}{1-N(\p)^{-1}}\Biggl(N(\p)^{-1}f|_v(1) \\
    &\quad +N(\mathfrak{d}_v)^{\frac{1}{2}}\int\limits_{|x_v|_v < 1}\dfrac{f|_v(x_v)-f|_v(1)}{|1-x_v|_v}\,dx_v \\ 
    &\quad +N(\mathfrak{d}_v)^{\frac{1}{2}}\int\limits_{|x_v|_v>1}\dfrac{f|_v(x_v)}{|1-x_v|_v}\,dx_v\Biggr)\\ 
 &= \log(N(\mathfrak{d}_v))f|_v(1)-\frac{\log(N(\p))}{1-N(\p)^{-1}}N(\mathfrak{d}_v)^{\frac{1}{2}}\int\limits_{|x_v|_v \neq 1}\dfrac{f|_v(x_v)}{|1-x_v|_v}\,dx_v \\
     &= \log(N(\mathfrak{d}_v))f(1)-\log(N(\p)) \sum_{n\neq 0}f(N(\p)^n)\min(1,N(\p)^n) \\
    &= \log(N(\mathfrak{d}_v))f(1)+ \frac{1}{2\pi i}\int_{\sigma-i\infty}^{\sigma+i\infty}\mathcal{M}^s(f)\,d\log\Biggl(\frac{1-N(\p)^{s-1}}{1-N(\p)^{-s}}\Biggr) \\
    &= \frac{1}{2\pi i}\int_{\sigma-i\infty}^{\sigma+i\infty}\mathcal{M}^s(f)\,d\log\Biggl(N(\mathfrak{d}_v)^{s-1/2}\frac{1-N(\p)^{s-1}}{1-N(\p)^{-s}} \Biggr).
\end{aligned}
\end{equation*}
\end{proof}

\subsection{Unramified case} Let, $K_v/\Qp$ be an unramified quadratic extension where $|\p|=N(\p)^{-1}$, with $N(\p)=p^2$ or $N(\p)=p$. For the sake of clarity and to make this article self-contained, we include the calculations here, unified for the unramified and inert cases, the calculation remains unchanged from the one described in \cite{Haran_1990} and \cite{Rodriguez_2008}. Specifically, the Riesz kernel is given by
\begin{equation*}
	R_v^s(x)=\dfrac{1-N(\p)^{-s}}{1-N(\p)^{s-1}}|x|_v^{s-1}\,dx_v \quad \text{Re } s>0, s\not\equiv 1 \pmod{\frac{2\pi i}{\log N(\p)}}.
\end{equation*}

As a distribution, it has a meromorphic continuation to all $s$, which is given by
\begin{multline*}
	R_v^s(\varphi)= \dfrac{1-N(\p)^{-1}}{1-N(\p)^{s-1}}\varphi(0) \\
	+\dfrac{1-N(\p)^{-s}}{1-N(\p)^{s-1}}\Bigl ( \int\limits_{|x_v|_v\leq 1} (\varphi(x)-\varphi(0))|x_v|_v^{s-1}\,dx_v +\int\limits_{|x_v|_v > 1} \varphi(x)|x_v|_v^{s-1}\,dx_v \Bigr).
\end{multline*}

In particular, for $\text{Re } s>0$:

\begin{align*}
	R_v^s(x)&= \dfrac{1-N(\p)^{-s}}{1-N(\p)^{s-1}}\int_{K_v}\varphi(x) |x|_v^{s-1}\,dx_v \\
	R_v^{-s}(x)&=\dfrac{1-N(\p)^{s}}{1-N(\p)^{-s-1}}\int_{K_v}\dfrac{\varphi(x)-\varphi(0)}{|x|_v^{1+s}}\,dx_v
\end{align*}
and for $s=0$, $R_{v}^0(\varphi)=\varphi(0)$, i.e. $R_{v}^0=\delta$.

Now, regard $R_v^s$ as an operator via convolution
\begin{multline*}
	R_v^s \ast \varphi(y)= \dfrac{1-N(\p)^{-1}}{1-N(\p)^{s-1}}\varphi(y) \\
	+\dfrac{1-N(\p)^{-s}}{1-N(\p)^{s-1}}\Bigl ( \int\limits_{|x_v|_v\leq 1} (\varphi(y+x)-\varphi(y))|x_v|_v^{s-1}\,dx_v +\int\limits_{|x_v|_v > 1} \varphi(y+x)|x_v|_v^{s-1}\,dx_v \Bigr)
\end{multline*}
once more, for $\text{Re } s>0$:
\begin{align*}
	R_v^s \ast \varphi(y)&= \dfrac{1-N(\p)^{-s}}{1-N(\p)^{s-1}}\int_{K_v}\varphi(y+x) |x|_v^{s-1}\,dx_v \\
	R_v^{-s} \ast \varphi(y)&=\dfrac{1-N(\p)^{s}}{1-N(\p)^{-s-1}}\int_{K_v}\dfrac{\varphi(y+x)-\varphi(y)}{|x|_v^{1+s}}\,dx_v
\end{align*}

If $\varphi$ is locally constant and has compact support, then $R_v^s \ast \varphi(y)$ remains locally constant, and its decay is characterized by $O(|x|_v^{\text{Re}(s)-1})$. This allows us to take the convolution $R_v^{s'} \ast (R_v^s \ast \varphi) = R_v^{s'+s} \ast \varphi$ under the condition $\text{Re }(s'+s)<1$ (see \cite{Kochubei}).
Specifically, when $\text{Re }s\geq 0$, the operator $R_v^{-s}$ forms a semigroup, and we compute next its infinitesimal generator below.
\begin{thm}
Let $f\in\mathscr{S}$ and $\sigma>0$, then
\begin{equation*}
		\left.\frac{\partial}{\partial s}\right|_{s=0} R_v^{-s}\ast f|_v(1)=\frac{1}{2\pi i}\int_{\sigma-i\infty}^{\sigma+i\infty}\mathcal{M}^s(f)\,d\log\Bigl(\frac{1-N(\p)^{s-1}}{1-N(\p)^{-s}} \Bigr)
	\end{equation*}
\end{thm}

\begin{proof}
\begin{equation*}
\begin{aligned}
	\left.\frac{\partial}{\partial s}\right|_{s=0} R_{v}^{-s}& * f|_v(1)= \left.\frac{\partial}{\partial s}\right|_{s=0} \Biggr\{\dfrac{1-N(\p)^{-1}}{1-N(\p)^{-s-1}}f|_v(1)+ \dfrac{1-N(\p)^{s}}{1-N(\p)^{-s-1}}\times\\
	&\quad \Bigl(\int\limits_{|x_v|_v\leq 1}\dfrac{f|_v(1-x_v)-f|_v(1)}{|x_v|_v^{1+s}}\,dx_v+\int\limits_{|x_v|_v>1}\dfrac{f|_v(1-x_v)}{|x_v|_v^{1+s}}\,dx_v \Bigr) \Biggr \}\\
	&=-\frac{\log(N(\p))}{1-N(\p)^{-1}}\Bigl(N(\p)^{-1}f|_v(1)\\
	&\quad +\int\limits_{|x_v|_v\leq 1}\dfrac{f|_v(1-x_v)-f|_v(1)}{|x_v|_v}\,dx_v+\int\limits_{|x_v|_v>1}\dfrac{f|_v(1-x_v)}{|x_v|_v}\,dx_v\Bigr)\\
	&=-\frac{\log(N(\p))}{1-N(\p)^{-1}}\Bigl(N(\p)^{-1}f|_v(1)\\
	&\quad +\int\limits_{|x_v|_v < 1}\dfrac{f|_v(x_v)-f|_v(1)}{|1-x_v|_v}\,dx_v+\int\limits_{|x_v|_v>1}\dfrac{f|_v(x_v)}{|1-x_v|_v}\,dx_v\Bigr)\\
	&=-\frac{\log(N(\p))}{1-N(\p)^{-1}}\int\limits_{|x_v|_v \neq 1}\dfrac{f|_v(x_v)}{|1-x_v|_v}\,dx_v\\
	&=-\log(N(\p)) \sum_{n\neq 0}f(N(\p)^n)\min(1,N(\p)^n)\\
	&=\frac{1}{2\pi i}\int_{\sigma-i\infty}^{\sigma+i\infty}\mathcal{M}^s(f)\,d\log\Bigl(\frac{1-N(\p)^{s-1}}{1-N(\p)^{-s}}\Bigr).
\end{aligned}
\end{equation*}	
\end{proof}

\subsection{Complex case}
Let $K_v=\CC$, in this case the Riesz kernel is given by 
\begin{equation*}
R_{\CC}^s(z)=\frac{\zeta_{\CC}(1-s)}{\zeta_{\CC}(s)}|z|_\infty ^{s-1}\,dz_\infty, \quad \text{Re } s>0, s\neq 1,2,\dotsc
\end{equation*}
As a distribution, it has meromorphic continuation to all $s$, given by
\begin{multline*}
	R_{\CC}^s(\varphi)=(2\pi)^{2s}\dfrac{\Gamma(1-s)}{\Gamma(1+s)}\varphi(0)+\\
	(2\pi)^{2s-1}\dfrac{\Gamma(1-s)}{\Gamma(s)}\Bigl(\int\limits_{|z|_\infty\leq 1}(\varphi(|z|_\infty)-\varphi(0))|z|_\infty ^{s-1}\,dz_\infty +\int\limits_{|z|_\infty >1}\varphi(|z|_\infty )|z|_\infty ^{s-1}\,dz_\infty \Bigr),
\end{multline*}

in polar coordinates
\begin{multline*}
	R_{\CC}^s(\varphi)=(2\pi)^{2s}\dfrac{\Gamma(1-s)}{\Gamma(1+s)}\varphi(0)+\\
	(2\pi)^{2s-1}\dfrac{\Gamma(1-s)}{\Gamma(s)}\Bigl(\int\limits_{r^2\leq 1}(\varphi(r^2)-\varphi(0))(r^2)^{s-1}\,2rdrd\theta+\int\limits_{r^2>1}\varphi(r^2)(r^2)^{s-1}\,2rdrd\theta\Bigr),
\end{multline*}
and making the change of variables $t=r^2$
\begin{multline*}
	R_{\CC}^s(\varphi)=(2\pi)^{2s}\dfrac{\Gamma(1-s)}{\Gamma(1+s)}\varphi(0)+\\
	+(2\pi)^{2s}\dfrac{\Gamma(1-s)}{\Gamma(s)}\Bigl(\int_0^1(\varphi(t)-\varphi(0))t^{s-1}\,dt+\int_1^{\infty}\varphi(t)t^{s-1}\,dt\Bigr),
\end{multline*}
in particular for $\text{Re}(s)>0$
\begin{align*}
	R_{\CC}^s(\varphi)&=(2\pi)^{2s}\dfrac{\Gamma(1-s)}{\Gamma(s)}\int_0^{\infty}\varphi(t)t^{s-1}\,dt\\
	R_{\CC}^{-s}(\varphi)&=(2\pi)^{-2s}\dfrac{\Gamma(1+s)}{\Gamma(-s)}\int_0^{\infty}\dfrac{\varphi(t)-\varphi(0)}{t^{s+1}}\,dt,
\end{align*}
and for $s=0$, $R_{\CC}^0(\varphi)=\varphi(0)$, i.e. $R_{\CC}^0=\delta$.

Now consider $R_{\CC}^s$ as an operator via convolution, we can form $R_{\CC}^s\ast \varphi$
\begin{multline*}
	R_{\mathbb{C}}^s \ast \varphi(y) = (2\pi)^{2s} \frac{\Gamma(1-s)}{\Gamma(1+s)}\varphi(|y|_\infty)\\
	+ (2\pi)^{2s-1}\frac{\Gamma(1-s)}{\Gamma(s)}\left( \int\limits_{|z|_\infty \leq 1}(\varphi(|y+z|_\infty)-\varphi(|y|_\infty))|z|_\infty^{s-1}\,dz_\infty \right. \\
	\left. + \int\limits_{|z|_\infty >1}\varphi(|y+z|_\infty)|z|_\infty^{s-1}\,dz_\infty \right).
\end{multline*}

Recall that the Riesz potentials in $\mathbb{R}^n$ are given by (see \cite[Chapter III]{stein}, \cite{Berg})
\begin{equation*}
	(I_{\alpha}f)(x)=\frac{\Gamma\left(\frac{n-\alpha}{2}\right)}{\pi^{n/2}2^{\alpha}\Gamma\left(\frac{\alpha}{2}\right)}\int_{\mathbb{R}^n} f(y)|x-y|^{\alpha-n}\,dy,
\end{equation*}
and, because $dz_\infty=2\,dx\,dy$, we have that $R_{\mathbb{C}}^s$ is the Riesz potential $I_{2s}$ in $\mathbb{R}^2$ multiplied by the factor $2^{2s}(2\pi)^{2s}$. Therefore, we have the semigroup property 
\begin{equation*}
	R_{\mathbb{C}}^{s'}\ast R_{\mathbb{C}}^s=R_{\mathbb{C}}^{s'+s} \quad \text{for } s,s'>0 \quad \text{and } s+s'<1,
\end{equation*}
and again, for $\text{Re}(s)\geq 0$, $R_{\mathbb{C}}^{-s}$ forms a holomorphic semigroup of operators with infinitesimal generator $\left.\frac{\partial}{\partial s}\right|_{s=0} R_{\mathbb{C}}^{-s}$. We will now proceed to compute this infinitesimal generator.

\begin{thm}
Let $f\in\mathscr{S}$ and $\sigma>0$, then
\begin{equation*}
	\left.\frac{\partial}{\partial s}\right|_{s=0} R_{\CC}^{-s} * f|_\infty(1)=
	\frac{1}{2\pi i}\int_{\sigma-i\infty}^{\sigma+i\infty}\mathcal{M}(f)(s)d\,\log\dfrac{\zeta_{\CC}(s)}{\zeta_{\CC}(1-s)}
\end{equation*}	
\end{thm}
\begin{proof}
\begin{equation*}
\begin{aligned}
	\left.\frac{\partial}{\partial s}\right|_{s=0} R_{\CC}^{-s}& * f|_\infty(1)= \left.\frac{\partial}{\partial s}\right|_{s=0} \Biggr\{ (2\pi)^{-2s}\dfrac{\Gamma(1+s)}{\Gamma(1-s)}f(1)\\
	&\quad +(2\pi)^{-2s-1}\dfrac{\Gamma(1+s)}{\Gamma(-s)}\Biggl(\,\int\limits_{r^2\leq 1}\dfrac{f(1+r^2)-f(1)}{(r^2)^{s+1}}\,2rdrd\theta\\
	&\hspace{4.6cm} +\int\limits_{r^2>1}\dfrac{f(1+r^2)}{(r^2)^{s+1}}\,2rdrd\theta \Biggr) \Biggr \}\\
	&=-2(\gamma+\log(2\pi))f(1)\\
	&-\dfrac{1}{2\pi}\Biggl\{ \,\int\limits_{r^2\leq 1}\dfrac{f(1+r^2)-f(1)}{r^2}\,2rdrd\theta
	+\int\limits_{r^2>1}\dfrac{f(1+r^2)}{r^2}\,2rdrd\theta \Biggr\}\\
	&\hspace{0.5cm}=-2(\gamma+\log(2\pi))f(1)\\
	&\hspace{0.5cm}-\int_{|1-r^2|<1}\dfrac{f(r^2)-f(1)}{|1-r^2|} \,2rdr
	-\int_{|1-r^2|>1}\dfrac{f(r^2)}{|1-r^2|}\,2rdr\\
	&=-2(\gamma+\log(2\pi))f(1)\\
	&\hspace{0.5cm}-\int_{|1-u|<1}\dfrac{f(u)-f(1)}{|1-u|} \,du
	-\int_{|1-u|>1}\dfrac{f(u)}{|1-u|}\,du\\
	&=-2(\gamma+\log(2\pi))f(1)-\int_1^{\infty}f(u)\,\dfrac{du}{u}-\int_1^{\infty}\dfrac{f(u)-f(1)}{u-1}\,\dfrac{du}{u}\\
	&\hspace{0.5cm} -\int_0^1f(u)\,du-\int_0^1\dfrac{uf(u)-f(1)}{1-u}\,du\\
	&=\frac{1}{2\pi i}\int_{(c)}\mathcal{M}(f)(s)\Biggl(-(\gamma+\log(2\pi))-\frac{1}{s}+\sum_{n\geq 1}\frac{1}{n}-\frac{1}{n+s}\\
	&\hspace{0.5cm} -(\gamma+\log(2\pi))-\frac{1}{1-s}+\sum_{n\geq 1}\frac{1}{n}-\frac{1}{n+1-s}\Biggr)\\
	&=\frac{1}{2\pi i}\int_{\sigma-i\infty}^{\sigma+i\infty}\mathcal{M}(f)(s)d\,\log\dfrac{\zeta_{\CC}(s)}{\zeta_{\CC}(1-s)}.
\end{aligned}
\end{equation*}		
\end{proof}

\section{Explicit Formula}\label{ExplicitFormula}
In this section, we will combine the Riesz kernels associated with each one of the places in the number field and demonstrate that, following Haran's approach for $\Q$, the Weil explicit formula in the case of a imaginary quadratic field can be interpreted in terms of the adelic Riesz potential.

\subsection{Global formula}
We shall rescale the infinitesimal generators $R_v^s$ and combine them to generate an operator defined on $\mathscr{S}$. Let $\phi_v^{\ast}(x)$ the characteristic function of $\mathcal{O}_v^{\times}$.

We define for $K_v$ unramified
\begin{equation*}
	c_v^s=R_v^{-s}\ast \phi_v^{\ast}(1)=1+O(s^2),
\end{equation*}
for $K_v$ totally ramified
\begin{equation*}
	c_v^s=N(\mathfrak{d}_v)^{-s} R_v^{-s}\ast \phi_v^{\ast}(1)=1+O(s^2),
\end{equation*}
and finally we set $c_\infty^s=1$

\begin{thm}
Let $f\in\mathscr{S}$, then
	\begin{equation*}
		\sum_{k \in K^{\times}} \left.\frac{\partial}{\partial s}\right|_{s=0} \Delta_{\mathbb{A}^{\times}}^s f(k)=\sum_{\zeta_K(s)=0} \mathcal{M}^s(f)-\mathcal{M}^0(f)-\mathcal{M}^1(f)
	\end{equation*}
	in this expression, the right-hand sum extends over all the zeros of $\zeta_K(s)$.
\end{thm}
\begin{proof}
It is sufficient to consider an elementary symmetric function $f=\otimes_{v}f_{v}$, thus
\begin{align*}
	\sum_{k \in K^{\times}} &\left.\frac{\partial}{\partial s}\right|_{s=0} \Delta_{\mathbb{A}^{\times}}^s f(k)=\sum_{k \in K^{\times}} \left.\frac{\partial}{\partial s}\right|_{s=0} \prod_{v}\Delta_{v}^s f_{v}(k)\\
	&=\sum_{k \in K^{\times}} \left.\frac{\partial}{\partial s}\right|_{s=0} \prod_{v}\Delta_{v}^s \bigl(f_{v}(kx_v)\bigr)(1)\\
	&=\sum_{k \in K^{\times}} \sum_v \left.\frac{\partial}{\partial s}\right|_{s=0} \prod_{v}\Delta_{v}^s \bigl(f_{v}(kx_{v})\bigr)(1)\cdot\prod_{v'\neq v}f_{v'}(k)\\
	&=\sum_{v} \left.\frac{\partial}{\partial s}\right|_{s=0} \Delta_{v}^s \bigl(\tilde{f}|_{v} \bigr)(1)\\
	&=\sum_{v} \left.\frac{\partial}{\partial s}\right|_{s=0} R_{v}^{-s}\ast \tilde{f}|_{v}(1) \\
	&=\sum_{v}\frac{1}{2\pi i}\int_{c-i\infty}^{c+i\infty} \mathcal{M}^s(\tilde{f}) d\,\log N(\mathfrak{d_v})^{s-\frac{1}{2}}\dfrac{\zeta_v(s)}{\zeta_v(1-s)}\\
	&=\sum_{v}\frac{1}{2\pi i}\int_{c-i\infty}^{c+i\infty} \mathcal{M}^s(\tilde{f}) d\,\log N(\mathfrak{d_v})^{s-\frac{1}{2}}\zeta_v(s)\\
	&\hspace{3cm}-\sum_{v}\frac{1}{2\pi i}\int_{c-i\infty}^{c+i\infty} \mathcal{M}^{1-s}(\tilde{f}) d\,\log\zeta_v(s)\quad \text{ for any $c>0$}\\
	&=\frac{1}{2\pi i}\int_{c-i\infty}^{c+i\infty} \mathcal{M}^s(\tilde{f}) d\,\log |D_K|^{s-\frac{1}{2}}\zeta_K(s)-\frac{1}{2\pi i}\int_{c-i\infty}^{c+i\infty} \mathcal{M}^{1-s}(\tilde{f}) d\,\log \zeta_K(s)\\
	&=\frac{1}{2\pi i}\int_{c-i\infty}^{c+i\infty} +\int_{1-c+i\infty}^{1-c-i\infty} \mathcal{M}^s(\tilde{f}) d\,\log |D_K|^{s-\frac{1}{2}}\zeta_K(s)\\
	&=\sum_{\zeta_K(s)=0} \mathcal{M}^s(f)-\mathcal{M}^0(f)-\mathcal{M}^1(f).
\end{align*}
\end{proof}

\section{Acknowledgements}
The authors wish to extend their sincere gratitude to UPTC (Universidad Pedagógica y Tecnológica de Colombia) and Universidad Nacional de Colombia for their invaluable contributions to for their invaluable contributions to this research.

\bibliographystyle{plain}
\bibliography{biblio}

\end{document}